\input amstex

\documentstyle{amsppt} \magnification=\magstep1

\TagsOnRight \NoBlackBoxes

\hoffset1 true pc \voffset2 true pc

\hsize36 true pc \vsize52 true pc

\tolerance=2000
\def\m1{^{-1}}
\def\ov1{\overline}

\catcode`\@=11
\def\logo@{}
\catcode`\@=\active

\topmatter

\title
Integral group ring of the first Mathieu simple group
\endtitle
\author
V.A.~Bovdi, A.B.~Konovalov
\endauthor
\abstract We investigated the classical Zassenhaus conjecture  for
the  normalized unit group of the integral group ring of  the
simple Mathieu group $M_{11}$. As a consequence, for this group we
confirm the  conjecture by Kimmerle about prime graphs.
\endabstract
\subjclass
Primary 16S34, 20C05, secondary 20D08
\endsubjclass
\thanks
The research  was supported by OTKA T 037202, T 038059 and
ADSI107(VUB).
\endthanks
\address
\hskip-\parindent {\rm  V.A.~Bovdi\newline Institute of
Mathematics, University of Debrecen\newline P.O.  Box 12, H-4010
Debrecen, Hungary\newline Institute of Mathematics and
Informatics, College of Ny\'\i regyh\'aza\newline S\'ost\'oi \'ut
31/b, H-4410 Ny\'\i regyh\'aza, Hungary
\newline E-mail: vbovdi\@math.klte.hu}
\newline
\bigskip
\hskip-\parindent {\rm A.B.~Konovalov
\newline Department of Mathematics, Vrije Universiteit Brussel
\newline Pleinlaan 2, B-1050 Brussel, Belgium
\newline E-mail: konovalov\@member.ams.org
\newline http://ukrgap.exponenta.ru/konoval.htm}
\endaddress
\endtopmatter

\document

\head
Introduction  and main results
\endhead

Let $V(\Bbb ZG)$ be the normalized unit group  of the integral
group ring $\Bbb ZG$ of a  finite group $G$. The following famous
conjecture was formulated by H.~Zassenhaus in \cite{15}:

\proclaim{(ZC)} Every torsion unit $u \in V(\Bbb ZG)$ is conjugated
within rational group algebra $\Bbb QG$ to an element of $G$.
\endproclaim

This conjecture is already confirmed for several classes of groups
but, in general, the problem remains open, and the counterexample is not known.

Various methods were developed to deal with this conjecture. One
of original ones was suggested by I.S.~ Luthar and I.B.S.~Passi
\cite{12, 13}, and further it was improved by M.Hertweck \cite{9}.
Using this method, the conjecture was proved for several new
classes of groups, in particular, for $S_5$ and
for some finite simple groups (see \cite{4, 9, 10, 12, 13}).

The Zassenhaus conjecture appeared to be very hard, and several weakened
variations of it were formulated (see, for example, \cite{3}). One of the most interesting modifications was suggested by W.~Kimmerle \cite{11}. Let us briefly introduce it now.

Let $G$ be a finite group. Denote by $\#(G)$ the set of all primes
dividing the order of $G$. Then the {\it Gruenberg-Kegel graph} (or the
{\it prime graph}) of G  is a graph $\pi(G)$ with vertices labelled by
primes from $\#(G)$, such that vertices $p$ and $q$ are adjacent
if and only if there is an element of order $pq$ in the group $G$.
Then the conjecture by Kimmerle can be formulated in the following way:

\proclaim{(KC)}
If $G$ is a finite group, then $\pi(G)=\pi(V(\Bbb ZG))$.
\endproclaim

For Frobenius groups and solvable groups this conjecture was confirmed in \cite{11}.
In the present paper we continue the investigation of {\bf (KC)}, and confirmed it
for the first simple Mathieu group $M_{11}$, using the Luthar-Passi method. Moreover,
this allows us to give a partial solution of {\bf (ZC)} for $M_{11}$.

Our main results are the following:

\proclaim {Theorem} Let $V(\Bbb ZG)$  be the normalized unit  group   of the integral group
ring $\Bbb ZG$, where $G$ is the
simple Mathieu group $M_{11}$. Let $u$ be  a torsion unit of $V(\Bbb ZG)$ of order $|u|$. Then:
\itemitem{(i)} if  $|u|\not=12$, then $|u|$ coincides  with the order of some element  $g\in G$;
\itemitem{(ii)} if $|u| \in \{2,3,5,11\}$, then $u$ is rationally
conjugated with some  $g\in G$;
\itemitem{(iii)} if $|u|=4$, then the tuple of partial augmentations
of $u$ belongs to the set
$$
\split \{\;  (\nu_{2a}, \nu_{3a}, \nu_{4a}, \nu_{6a},\;& \nu_{5a},
\nu_{8a}, \nu_{8b}, \nu_{11a}, \nu_{11b}) \in \Bbb Z^9
\quad \mid \nu_{kx}=0, \\
&kx\not\in \{2a, 4a\},  \quad (\nu_{2a}, \nu_{4a} ) \in \{\; ( 0,
1 ), \; ( 2, -1 ) \; \} \; \};
\endsplit
$$
\itemitem{(iv)} if $|u|=6$,  then  the tuple of partial augmentations
of $u$ belongs to the set
$$
\split \{\; (\nu_{2a}, \nu_{3a}, \nu_{4a}, \nu_{6a}, \nu_{5a},
\nu_{8a}, \nu_{8b}, \nu_{11a}, \nu_{11b})\in \;&\Bbb Z^9
\quad \mid \nu_{kx}=0,\\
kx\not\in \{2a, 3a, 6a\},\quad  (\nu_{2a}, \nu_{3a}, \nu_{6a}) \in
\{\;
&( -2, 3, 0 ),  ( 0, 0,1),\\
( 0, 3, -2 ),\;  &( 2, -3, 2 ),\; ( 2, 0, -1) \; \} \quad  \};
\endsplit
$$
\itemitem{(v)} if $|u|=8$, then the tuple of partial augmentations
of $u$ belongs to the set
$$
\split \{\; (\nu_{2a},  \nu_{3a}, \nu_{4a}, \nu_{6a}, \nu_{5a},
\nu_{8a},\nu_{8b}, \nu_{11a}, \nu_{11b}) \in &\Bbb Z^9
\quad \mid \nu_{kx}=0,\\
kx\not\in \{4a, 8a, 8b\},\quad (\nu_{4a}, \nu_{8a}, \nu_{8b} ) \in
\{\; &(0,0,1 ),\; (0,1,0),\\
&(2,-1,0),\;  (2,0,-1) \; \}\quad \};
\endsplit
$$
\itemitem{(vi)} if $|u|=12$,  then the tuple of partial
augmentations of $u$ can not belong to the set
$$
\split \Bbb Z^9 \setminus \{\; (\nu_{2a},\;&  \nu_{3a}, \nu_{4a},
\nu_{6a}, \nu_{5a}, \nu_{8a}, \nu_{8b}, \nu_{11a},
\nu_{11b})\in\Bbb Z^9 \quad \mid\quad \nu_{kx}=0,\\
& kx\not\in \{2a, 4a, 6a\},\quad  (\nu_{2a}, \nu_{4a}, \nu_{6a})
\in \{\; (-1,1,1),\; (1,1,-1) \; \}\quad  \}.
\endsplit
$$
\endproclaim

\proclaim {Corollary} Let $V(\Bbb ZG)$ be the normalized unit  group  of the integral
group ring $\Bbb ZG$, where $G$ is the simple Mathieu group $M_{11}$.  Then
$\pi(G)=\pi(V(\Bbb ZG))$, where  $\pi(G)$ and $\pi(V(\Bbb ZG))$ are prime graphs of
$G$ and $V(\Bbb ZG)$, respectively. Thus,  for $M_{11}$ the conjecture by Kimmerle
is true.
\endproclaim

\head
Notation and preliminaries
\endhead

Let $u=\sum\alpha_gg$ be a normalized torsion unit of order $k$
and let $\nu_i=\varepsilon_{C_i}(u)$ be a partial augmentation of
$u$. By  S.D.~Berman's Theorem \cite{2} we have that $\nu_1$=0 and
$$
\nu_2+\nu_3+\cdots+\nu_m=1.\tag1
$$
For any character $\chi$ of $G$ of degree $n$, we got that
$\chi(u)=\sum_{i=2}^m\nu_i\chi(h_i)$,
where $h_i$ is a representative  of a conjugacy class $ C_i$.

We need the following results:

\proclaim {Proposition 1} (see \cite{12})
Suppose that $u$ is an element of
$\Bbb ZG$ of order $k$. Let $z$ be a primitive $k$-th root of
unity. Then for every integer $l$ and any character $\chi$ of $G$,
the number
$$
\mu_l(u,\chi)=\frac{1}{k} \sum_{d|k}Tr_{\Bbb Q(z^d)/\Bbb Q}
\{\chi(u^d)z^{-dl}\}\tag2
$$
is a non-negative integer.
\endproclaim

\proclaim {Proposition 2}  (see  \cite{6})
Let $u$ be a torsion unit in  $V(\Bbb ZG)$.
Then the order of $u$ divides the exponent of $G$.
\endproclaim

\proclaim{Proposition 3}(see \cite{12} and Theorem 2.7 in \cite{14})
Let $u$ be a torsion unit of $V(\Bbb ZG)$.
Let $C$ be a conjugacy class of $G$. If $p$ is a prime
dividing  the order of a representative of $C$ but not the
order of $u$ then the partial augmentation $\varepsilon_C(u)=0$.
\endproclaim
M.~Hertweck (see \cite{10}, Proposition 3.1; \cite{9}, Lemma
5.6) obtained the next result:

\proclaim {Proposition 4} Let $G$ be a finite group and let $u$ be
a torsion unit in $V(\Bbb ZG)$.
\itemitem{(i)}
If $u$ has order $p^n$, then $\varepsilon_x(u)=0$
for every $x$ of $G$ whose $p$-part is of order strictly
greater than $p^n$.
\itemitem{(ii)}
If $x$ is an element of $G$ whose $p$-part, for some
prime $p$, has order strictly greater than the order of the $p$-part
of $u$, then $\varepsilon_x(u)=0$.
\endproclaim
Note that the first part of the Proposition 4 gives a partial
answer to the conjecture by A.~Bovdi (see \cite{1}). Also
M.~Hertweck (\cite{9}, Lemma 5.5) gives a complete answer to the
same conjecture in the case when $ G = \text{\rm{PSL}}(2,\Bbb F)$,
where $ \Bbb F = \text{\rm{GF}}(p^k)$.

In the sequel of the paper  for the partial augmentation  $\nu_{i}$
we shall also use the notation $\nu_{kx}$, where $k$ is the order of the
representative of the $i$-th conjugacy class,
and $x$ is a distinguishing letter for  this particular class
with elements of order $k$.

\proclaim{Proposition 5}(see \cite{12} and Theorem 2.5 in \cite{14})
Let $u$ be a torsion unit of $V(\Bbb ZG)$ of order $k$. Then $u$ is conjugate
in $\Bbb QG$ to an element $g \in G$ if and only if for each $d$
dividing $k$ there is precisely one conjugacy class $C_{i_d}$ with
partial augmentation $\varepsilon_{C_{i_d}}(u^d)  \neq 0 $.
\endproclaim

\proclaim{Proposition 6}(see \cite{6})
Let $p$ be a prime, and let $u$ be a torsion unit of $V(\Bbb ZG)$ of order $p^n$.
Then for $m \ne n$ the sum of all partial augmentations of $u$ with respect to
conjugacy classes of elements of order $p^m$ is divisible by $p$.
\endproclaim

The Brauer character table modulo $p$ of the group $M_{11}$
will be denoted by $\frak{BCT}{(p)}$.

\head
Proof of the Theorem
\endhead

It is well known \cite{7,8}
that $|G|=2^4 \cdot 3^2 \cdot 5 \cdot 11$ and  $exp(G) =1320$.
The  character table of $G$, as well as the Brauer character tables
$\frak{BCT}{(p)}$, where $p \in \{2,3,5,11\}$, can be found using
the computational algebra system GAP \cite{7}.

Since the group $G$ possesses elements of orders
$2$, $3$, $4$, $5$, $6$, $8$ and $11$,
first of all we shall investigate units of these orders.
After this, by Proposition 2, the order of each torsion unit divides the
exponent of $G$, so it will be enough to consider units of orders
$10$, $12$, $15$, $22$, $24$, $33$ and $55$, because if $u$ will be a
unit of another possible order, then there is $t \in \Bbb N$ such that
$u^t$ has an order from this list. We shall prove that units of all these
orders except 12 do not appear in $V(\Bbb ZG)$. For units of order 12 we are
not able to prove this, but we reduce this question to only two cases.

Let $u \in V(\Bbb ZG)$ has the order $k$.
By S.D.~Berman's Theorem \cite{2} and Proposition 3 we have $\nu_{1a}=0$ and
$$
\matrix
\nu_{3a}=\nu_{5a}=\nu_{6a}=\nu_{11a}=\nu_{11b}=0        & \text{when} \quad  k=2,4,8;\\
\nu_{2a}=\nu_{4a}=\nu_{5a}=\nu_{6a}=\nu_{8a}=\nu_{8b}=\nu_{11a}=\nu_{11b}=0
                                                        & \text{when} \quad  k=3;\\
\nu_{2a}=\nu_{3a}=\nu_{4a}=\nu_{6a}=\nu_{8a}=\nu_{8b}=\nu_{11a}=\nu_{11b}=0
                                                        & \text{when} \quad  k=5;\\
\nu_{5a}=\nu_{11a}=\nu_{11b}=0                          & \text{when} \quad  k=6;\\
\nu_{2a}=\nu_{3a}=\nu_{4a}=\nu_{5a}=\nu_{6a}=           & \\
\nu_{8a}=\nu_{8b}=0                                     & \text{when} \quad  k=11;\\
\nu_{3a}=\nu_{6a}=\nu_{11a}=\nu_{11b}=0                 & \text{when} \quad  k=10;\\
\nu_{5a}=\nu_{11a}=\nu_{11b}=0                          & \text{when} \quad  k=12;\\
\nu_{2a}=\nu_{4a}=\nu_{6a}=\nu_{8a}=\nu_{8b}=\nu_{11a}=\nu_{11b}=0
                                                        & \text{when} \quad  k=15;\\
\nu_{2a}=\nu_{3a}=\nu_{4a}=\nu_{{5a}}=\nu_{6a}=0        & \text{when} \quad  k=22;\\
\nu_{5a}=\nu_{11a}=\nu_{11b}=0                          & \text{when} \quad  k=24;\\
\nu_{2a}=\nu_{4a}=\nu_{{5a}}=\nu_{6a}=\nu_{8a}=\nu_{8b}=0
                                                        & \text{when} \quad  k=33;\\
\nu_{2a}=\nu_{3a}=\nu_{4a}=\nu_{6a}=\nu_{8a}=\nu_{8b}=0 & \text{when} \quad  k=55.\\
\endmatrix
\tag3
$$
It follows immediately that by Proposition 5 units of orders 3 and
5 are rationally conjugated with some element of $G$.

Now we consider each case separately:

\noindent $\bullet$ Let $u$ be an involution. Then using
(3) and Proposition 4 we obtain that
$\nu_{4a}=\nu_{8a}=\nu_{8b}=0$, so $\nu_{2a}=1$.

\noindent $\bullet$ Let $u$ be a unit of order $4$. Then by (3)
and Proposition 4 we have $\nu_{2a}+\nu_{4a}=1$. By (2),\quad
$\mu_0(u,\chi_3) = \textstyle\frac{1}{4} ( -4 \nu_{2a} + 8 ) \ge
0$\quad  and \quad  $\mu_2(u,\chi_3) = \textstyle\frac{1}{4} ( 4
\nu_{2a} + 8 ) \ge 0 $, so $\nu_{2a} \in \{ -2, -1, 0, 1, 2\}$.
Now using inequalities
$$
\mu_0(u,\chi_5) = \textstyle\frac{1}{4} ( 6 \nu_{2a} -2 \nu_{4a} +
14 ) \ge 0 ;\quad \mu_2(u,\chi_5) = \textstyle\frac{1}{4} ( -6
\nu_{2a} + 2 \nu_{4a} + 14 ) \ge 0,
$$
we get that there are only two integral solutions
$(\nu_{2a},\nu_{4a}) \in \{\; (0,1),\; (2,-1)\; \}$ satisfying (1)
and Proposition 6, such that all $\mu_i(u,\chi_j)$ are
non-negative integers.

\noindent$\bullet$ Let $u$ be a unit of order $6$. Then by (1),
(2) and Proposition 4 we get
$$
\nu_{2a}+\nu_{3a}+\nu_{6a}=1.
$$
Now, using $\frak{BCT}{(11)}$ from the system of inequalities
$\mu_0(u,\chi_6) = \textstyle\frac{1}{6}(-4 \nu_{3a}+12) \geq 0$
and $\mu_3(u,\chi_6) = \textstyle\frac{1}{6}( 4 \nu_{3a}+12) \geq
0$, we have that $\nu_{3a} \in \{-3,0,3\}$. Furthermore, from the
system of inequalities
$$
\split
\mu_3(u,\chi_2) &= \textstyle\frac{1}{6}( -2 \nu_{2a} + 4 \nu_{6a} + 8 ) \geq 0;\\
\mu_0(u,\chi_2) &= \textstyle\frac{1}{6}( 2  \nu_{2a} - 4 \nu_{6a} + 10 ) \geq 0;\\
\mu_1(u,\chi_2) &= \textstyle\frac{1}{6}(    \nu_{2a} - 2 \nu_{6a} + 8 ) \geq 0,\\
\endsplit
$$
we get that $\nu_{2a}-2\nu_{6a}\in\{-2,4\}$,
so $\nu_{6a} \in \{-2,-1,0,1,2\}$.
Using inequalities
$$
\split
\mu_0(u,\chi_3) &= \textstyle\frac{1}{6}(-4\nu_{2a}+2\nu_{3a}+2\nu_{6a}+10) \geq 0;\\
\mu_2(u,\chi_3) &= \textstyle\frac{1}{6}( 2\nu_{2a}-\nu_{3a}-\nu_{6a}+7) \geq 0,\\
\endsplit
$$
we obtain only the following integral solutions
$(\nu_{2a}, \nu_{3a}, \nu_{6a})$:
$$
 \{\; \quad ( -2, 3, 0 ),\quad  ( 0, 0, 1 ), \quad ( 0, 3, -2 ),
\quad  ( 2, -3, 2 ),\quad  ( 2, 0, -1 ) \quad  \; \},
\tag4
$$
such that all $\mu_i(u,\chi_j)$ are non-negative integers.

Using the GAP package LAGUNA \cite{5}, we tested all possible
$\mu_i(u,\chi_j)$ for all tuples $(\nu_{2a},\nu_{3a},\nu_{6a})$
from (4), and were not able to produce a contradiction. Thus, in
this case, as well as in the case of elements of order 4, the
Luthar-Passi method is not enough to prove the rational conjugacy.

\noindent $\bullet$ Let $u$ be a unit of order $8$. By (3) and
Proposition 4 we have
$$
\nu_{2a}+\nu_{4a}+\nu_{8a}+\nu_{8b}=1.
$$
Since $|u^2|=4$, by (4) it yields that
$$
\chi_j(u^2)=\overline\nu_{2a}\chi_j{(2a)}+\overline\nu_{4a}\chi_j{(4a)}.
$$
Now using $\frak{BCT}{(3)}$ by (2) in the case when
$(\overline\nu_{2a},\overline\nu_{4a})=(0,1)$ we get
$$
\split
\mu_0(u,\chi_5) &= \textstyle\frac{1}{8} (-8\nu_{2a} + 8) \geq 0; \qquad\qquad
\mu_4(u,\chi_5) = \textstyle\frac{1}{8} ( 8\nu_{2a} +  8) \geq 0,\\
\mu_0(u,\chi_4) &= \textstyle\frac{1}{8} ( 8\nu_{2a} + 8\nu_{4a}
+ 16 ) \geq 0; \quad
\mu_4(u,\chi_4) = \textstyle\frac{1}{8} ( -8\nu_{2a}-8\nu_{4a}  + 16 ) \geq 0; \\
\mu_1(u,\chi_2) &= \textstyle\frac{1}{8} (4\nu_{8a}-4\nu_{8b} + 4 ) \geq 0; \qquad
\mu_4(u,\chi_7) = \textstyle\frac{1}{8} (-8\nu_{8a}-8\nu_{8b} + 24 ) \geq 0;\\
\mu_5(u,\chi_2) &= \textstyle\frac{1}{8} (-4\nu_{8a}+ 4\nu_{8b} + 4 ) \geq 0; \quad
\mu_0(u,\chi_7) = \textstyle\frac{1}{8} ( 8\nu_{8a}+ 8\nu_{8b} + 24 ) \geq 0. \\
\endsplit
$$
It follows that
\quad  $-1\leq \nu_{2a}\leq 1$,
\quad  $-3\leq \nu_{4a}\leq 3$,
\quad  $-2\leq \nu_{8a},\nu_{8b}\leq 2$.\quad
Considering additional inequality
$$
\mu_0(u,\chi_2) = \textstyle\frac{1}{8}
(4 \nu_{2a} - 4 \nu_{4a} - 4 \nu_{8a} - 4 \nu_{8b} + 4) \geq 0,
$$
and using Proposition 6, it is easy to check that this system
has the following integral solutions
$(\nu_{2a}, \nu_{4a}, \nu_{8a},\nu_{8b})$:
$$
\{ \; ( 0, 2, 0, -1 ),\; ( 0, 2, -1, 0 ),\;  ( 0, -2, 1, 2 ),\;
      ( 0, -2, 2, 1 ),\;  ( 0, 0, 1, 0 ),\;  ( 0, 0, 0, 1 )  \;
      \}.\tag5
$$
In the case when $(\overline\nu_{2a},\overline\nu_{4a})=(2,-1)$
using $\frak{BCT}{(3)}$ by (2) we obtained that
$\mu_0(u,\chi_5)=-\mu_4(u,\chi_5) =-\nu_{2a}=0$\quad and
$$
\split
\mu_0(u,\chi_4) &= \textstyle\frac{1}{8} (8\nu_{4a}  + 16 ) \geq 0; \qquad\quad
\mu_4(u,\chi_4)  = \textstyle\frac{1}{8} (-8\nu_{4a}  + 16 ) \geq 0; \\
\mu_1(u,\chi_2) &= \textstyle\frac{1}{8} (4\nu_{8a}-4\nu_{8b} + 4 ) \geq 0; \quad
\mu_5(u,\chi_2)  = \textstyle\frac{1}{8} (-4\nu_{8a}+ 4\nu_{8b} + 4) \geq 0. \\
\endsplit
$$
It is easy to check that
\; $-2\leq \nu_{4a}, \nu_{8a},\nu_{8b}\leq 2$, \;
and this system has the following integral solutions:
$$
\{\; ( 0, 2, -1, 0 ),\; ( 0, 2, 0, -1 ),\;  ( 0, 0, 0, 1 ),\; ( 0,
-2, 2, 1 ),\;  ( 0, 0, 1, 0 ),\;  ( 0, -2, 1, 2 )\; \}.\tag6
$$
Now using $\frak{BCT}{(11)}$ in the case when
$(\overline\nu_{2a},\overline\nu_{4a})=(0,1)$,
by (2) we get
$$
\split
\mu_0(u,\chi_3) &= \textstyle\frac{1}{8} (-8\nu_{2a} + 8 ) \geq 0; \quad
\mu_4(u,\chi_3)  = \textstyle\frac{1}{8} (8\nu_{2a} + 8 ) \geq 0;\\
\mu_1(u,\chi_3) &= \textstyle\frac{1}{8} (4t + 12) \geq 0; \qquad
\mu_5(u,\chi_3)  = \textstyle\frac{1}{8} (-4t + 12 ) \geq 0; \\
\mu_0(u,\chi_2) &= \textstyle\frac{1}{8} (4v-4w + 12 ) \geq 0;
\quad
\mu_4(u,\chi_2) = \textstyle\frac{1}{8} (-4v+ 4w+ 12 ) \geq 0; \\
\mu_0(u,\chi_5) &= \textstyle\frac{1}{8} (4z-4w+ 12 ) \geq 0;
\quad
\mu_4(u,\chi_5) = \textstyle\frac{1}{8} (-4z+ 4 w + 12 ) \geq 0,\\
\endsplit
$$
where $t=\nu_{8a}-\nu_{8b}$, $z=3\nu_{2a}-\nu_{4a}$,
$v=\nu_{2a}+\nu_{4a}$ and $w=\nu_{8a}+\nu_{8b}$.
From this follows that \; $-1\leq \nu_{2a} \leq 1$, \; \; $-8\leq
\nu_{2a} \leq 10$, \; \; $-4\leq \nu_{8a},\nu_{8b}\leq 4$, \; and
using Proposition 6, it is easy to check that this system has the
following integral solutions $(\nu_{2a}, \nu_{4a},
\nu_{8a},\nu_{8b})$:
$$
\split \{\quad
& ( 0, 2, 1, -2 ), \; ( 0, 2, -2, 1 ), \; ( 0, 2, 0, -1 ), \; ( 0, 2, -1, 0 ), \; \\
& ( 0, 0, 1, 0 ), \; ( 0, 0, 2, -1 ), \; ( 0, 0, 0, 1 ), \; ( 0,
0, -1, 2 ) \quad  \}.
\endsplit
\tag7
$$
In the case when $(\overline\nu_{2a},\overline\nu_{4a})=(2,-1)$
first using using
$\frak{BCT}{(11)}$ by (2) we get $\mu_0(u,\chi_3)=-\mu_4(u,\chi_3)
=-\nu_{2a}=0$ \quad and
$$
\split
\mu_1(u,\chi_3) &= \textstyle\frac{1}{8} (4\nu_{8a}-4\nu_{8b}+ 12 ) \geq 0; \\
\mu_5(u,\chi_3) &= \textstyle\frac{1}{8} (-4\nu_{8a}+4\nu_{8b}+ 12 ) \geq 0; \\
\mu_0(u,\chi_2) &= \textstyle\frac{1}{8} ( 4\nu_{4a}-4\nu_{8a}-4\nu_{8b}+ 12 ) \geq 0; \\
\mu_4(u,\chi_2) &= \textstyle\frac{1}{8} ( -4\nu_{4a}+ 4\nu_{8a}+ 4\nu_{8a}+ 12 ) \geq 0; \\
\mu_0(u,\chi_5) &= \textstyle\frac{1}{8} ( -4\nu_{4a}-4\nu_{8a}-4\nu_{8a}+ 28 ) \geq 0; \\
\mu_4(u,\chi_5) &= \textstyle\frac{1}{8} ( 4\nu_{4a}+ 4\nu_{8a}+
4\nu_{8a}+ 28 ) \geq 0.
\endsplit
$$
It is easy to check that \; $-7 \leq \nu_{4a} \leq 9$, \; \;
$-4\leq \nu_{8a},\nu_{8b}\leq 4$, and the system has the following
integral solutions $(\nu_{2a}, \nu_{4a}, \nu_{8a},\nu_{8b})$:
$$
\split \{\quad
& ( 0, 2, -1, 0 ), \; ( 0, 0, -1, 2 ), \; ( 0, 2, 0, -1 ), \; ( 0, 0, 0, 1 ), \\
& ( 0, 0, 2, -1 ), \; ( 0, 2, -2, 1 ), \; ( 0, 0, 1, 0 ), \; ( 0,
2, 1, -2 ) \quad  \}.
\endsplit
\tag8
$$
It follows from (5)-(8) that only four solutions appear in
both cases when $p=3$ and $p=11$ are the following ones:
$\nu_{2a}=0$ and
$$
(\nu_{4a},\nu_{8a},\nu_{8b}) \in \{\quad ( 0, 0, 1 ), ( 0, 1, 0 ),
( 2, -1, 0 ), ( 2, 0, -1 ) \quad \}.
$$
Again, we were not able to produce a contradiction computing all
possible $\mu_i(u,\chi_j)$ for all listed above tuples
$(\nu_{4a},\nu_{8a},\nu_{8b})$ for the ordinary character table of
$G$ as well as for $\frak{BCT}{(p)}$, where $p \in \{ 3,5,11 \}$.

\noindent $\bullet$ Let $u$ be a unit of order $11$. Then using
$\frak{BCT}{(3)}$  by (2) we have
$$
\split
\mu_1(u,\chi_2) & = \textstyle\frac{1}{11} ( 6 \nu_{11a}-5  \nu_{11b} + 5 )\geq 0;\\
\mu_2(u,\chi_2) & = \textstyle\frac{1}{11} ( -5 \nu_{11a}+ 6 \nu_{11b} + 5 )\geq 0,\\
\endsplit
$$
which has  only the following  trivial solutions
$(\nu_{11a},\nu_{11b}) = \{\;  (1,0),\;(0,1)\;  \}$.

For all above mentioned cases except elements of orders 4, 6 and 8
we obtained that there is precisely one conjugacy class with
non-zero partial augmentation. Thus, by Proposition 5, part
(ii) of the Theorem is proved.

It remains to prove parts (i) and (vi), considering units
of $V(\Bbb ZG)$ of orders $10$, $12$, $15$, $22$, $24$, $33$ and $55$.

\noindent $\bullet$ Let $u$ be a unit of order $10$. Then by (1),
(3) and Proposition 4 we get $\nu_{2a} + \nu_{5a} = 1$. Using
$\frak{BCT}{(3)}$ by (2) we have  the system of inequalities
$$
\split
\mu_5(u,\chi_4) & = \textstyle\frac{1}{10} ( -8 \nu_{2a} + 8 ) \geq 0; \quad \\
\mu_0(u,\chi_4) & = \textstyle\frac{1}{10} (  8  \nu_{2a} + 12 ) \geq  0; \\
\mu_2(u,\chi_2) & = \textstyle\frac{1}{10} ( -\nu_{2a} + 6 ) \geq  0,
\endsplit
$$
that has no integral solutions such that $\mu_5(u,\chi_4), \;
\mu_0(u,\chi_4), \; \mu_2(u,\chi_2) \in \Bbb Z$.

\noindent $\bullet$ Let $u$ be a unit of order $12$. By (1), (3)
and Proposition 4, we obtain that
$$
\nu_{2a}+ \nu_{3a} + \nu_{4a} +\nu_{6a} = 1.
$$
Since $|u^2|=6$ and $|u^3|=4$, by (4)  it yields that
$$
\chi_j(u^2)=\overline\nu_{2a}\chi_j{(2a)}+\overline\nu_{3a}\chi_j{(3a)}+\overline\nu_{6a}\chi_j{(6a)}
$$
and\quad
$\chi_j(u^3)=\widetilde\nu_{2a}\chi_j{(2a)}+\widetilde\nu_{4a}\chi_j{(4a)}$.\quad

Consider the following four  cases from parts (ii) and
(iii) of the Theorem:

\noindent {\bf 1.} Let $(\overline\nu_{2a}, \overline\nu_{3a},
\overline\nu_{6a}) \in \{\; ( 0, 0, 1 ),\; ( 2, -3, 2 ),\; ( -2,
3, 0 )\; \}$ and suppose that
$(\widetilde\nu_{2a},\widetilde\nu_{4a}) \in \{\; (0,1),\;
(2,-1)\; \}$. Then by (2) we have $\mu_0(u,
\chi_2)=\textstyle\frac{1}{2}\not\in \Bbb Z$, a contradiction.

\noindent {\bf 2.}  Let $(\overline\nu_{2a}, \overline\nu_{3a},
\overline\nu_{6a})=( 2, 0, -1 )$ and
$(\widetilde\nu_{2a},\widetilde\nu_{4a}) \in \{\; (0,1),(2,-1)\;
\} $. Then by (2) we have $\mu_1(u,
\chi_3)=\textstyle\frac{1}{2}\not\in \Bbb Z$, a contradiction.

\noindent {\bf 3.}  Let $(\overline\nu_{2a}, \overline\nu_{3a},
\overline\nu_{6a})=(0,3,-2)$ and
$(\widetilde\nu_{2a},\widetilde\nu_{4a})=(0,1)$.  According to
(2), $\mu_6(u,\chi_6)=-\mu_0(u,\chi_6) = \frac{2}{3}\nu_{3a}=0$,
\quad so $\nu_{3a}=0$, and we have the system
$$
\split \mu_2(u, \chi_5)&=\textstyle\frac{1}{12}( 6 \nu_{2a}-2
\nu_{4a}+8) \geq 0; \qquad \quad
\mu_4(u, \chi_5)=\textstyle\frac{1}{12}(-6 \nu_{2a}+2 \nu_{4a}+4) \geq 0; \\
\mu_2(u, \chi_3)&=\textstyle\frac{1}{12}(-4 \nu_{2a}+2 \nu_{6a}+6)
\geq 0; \qquad
\mu_4(u, \chi_3)=\textstyle\frac{1}{12}( 4 \nu_{2a}-2 \nu_{6a}+6) \geq 0; \\
\mu_4(u, \chi_2)&=\textstyle\frac{1}{12}(-4 \nu_{2a}-4 \nu_{4a} +
2 \nu_{6a} +10) \geq 0; \\
\mu_2(u, \chi_2)&=\textstyle\frac{1}{12}(4 \nu_{2a}+4 \nu_{4a} - 2 \nu_{6a} +2) \geq 0, \\
\endsplit
$$
that has only two solutions $\{\; (-1,0,1,1),\;  (1,0,1,-1)\;
\}$ with $\mu_i(u, \chi_j) \in \Bbb Z$.

\noindent {\bf 4.}  Let $(\overline\nu_{2a}, \overline\nu_{3a},
\overline\nu_{6a})=(0,3,-2)$ and
$(\widetilde\nu_{2a},\widetilde\nu_{4a})=(2,-1)$. Using (2), we
get $\mu_6(u,\chi_6)=-\mu_0(u,\chi_6) = \frac{2}{3}\nu_{3a}=0$, \;
so $\nu_{3a}=0$. Put $t=2\nu_{2a}-\nu_{6a}$, then by (2)
$$
\mu_0(u,\chi_3) = \textstyle\frac{1}{12} (-4t + 4) \geq 0; \qquad
\mu_4(u,\chi_3) = \textstyle\frac{1}{12} (2t-2) \geq 0,
$$
so $2\nu_{2a}-\nu_{6a}=1$. Now by (2) we have
$$
\mu_2(u,\chi_5)=\textstyle\frac{1}{12} (2(3\nu_{2a} -\nu_{4a})-8 )
\geq 0, \quad \mu_0(u,\chi_9) = \textstyle\frac{1}{12} (
-4(3\nu_{2a}-\nu_{4a}) + 28) \geq 0,
$$ and
$3\nu_{2a}-\nu_{4a}=4$. Using (1) we obtain that
$\nu_{2a}=\nu_{4a}=-\nu_{6a}=1$. Finally,
$\mu_4(u,\chi_9)=\textstyle\frac{1}{12} (6\nu_{2a}-2\nu_{4a}+
28)=\textstyle\frac{8}{3}\not\in \Bbb Z$, a contradiction. Thus,
Part (vi) of the Theorem is proved.

\noindent $\bullet$ Let $u$ be a unit of order $15$. Then by (1)
and (3) we have $\nu_{3a} +\nu_{5a}=1$. Now using the character
table of $G$, by (2) we get the system of inequalities
$$
\mu_0(u,\chi_2) = \textstyle\frac{1}{15} ( 8\nu_{3a} + 12 )\geq 0;
\quad \mu_5(u,\chi_2) =  \textstyle\frac{1}{15} ( -4\nu_{3a} + 9
)\geq 0,
$$
that has no integral solutions such that $\mu_0(u,\chi_2), \;
\mu_5(u,\chi_2) \in \Bbb Z$.

\noindent $\bullet$  Let $u$ be a unit of order $22$. Then by (1),
(3) and Proposition 4 we obtain that
$$
\nu_{2a}+\nu_{11a}+\nu_{11b}=1.
$$
In (2) we need to consider two
cases: $\chi(u^2)=\chi{(11a)}$ and $\chi(u^2)=\chi{(11b)}$, but in
both cases by (2) we have
$$
\mu_0(u,\chi_2) = -\mu_{11}(u,\chi_2)= \textstyle\frac{1}{22}(20
\nu_{2a} -10 \nu_{11a} -10\nu_{11b} + 2)=0.
$$
It yields $10 \nu_{2a} - 5 \nu_{11a} - 5 \nu_{11b} = -1$, that has
no integral solutions.

\noindent $\bullet$ Let $u$ be a unit  of order $24$. Then by (1)
and (3) we have
$$
\nu_{2a}+\nu_{3a}+\nu_{4a}+\nu_{6a}+\nu_{8a}+\nu_{8b}=1.
$$
Since $|u^2|=12$, $|u^4|=6$, $|u^3|=8$, $|u^6|=4$, and $G$ has two
conjugacy classes of elements of order $8$, we need to consider 40
various cases defined by Parts (iii)--(vi) of the
Theorem.

Let $(\overline\nu_{2a}, \overline\nu_{3a}, \overline\nu_{6a}) \in
\{\;  ( 0, 3, -2 ),\; ( 2, 0, -1)\; \}$, where
$$
\chi_j(u^4)=\overline\nu_{2a}\chi_j{(2a)}+\overline\nu_{3a}\chi_j{(3a)}+\overline\nu_{6a}\chi_j{(6a)}.
$$
According to  (2) we have that $\mu_1(u,
\chi_2)=\textstyle\frac{1}{2}\not\in\Bbb Z$, a contradiction. In
the remaining cases similarly we obtain that $\mu_1(u,
\chi_2)=\textstyle\frac{1}{4}\not\in\Bbb Z$, that is also a
contradiction.

\noindent $\bullet$ Let $u$ be a unit  of order $33$.  Then by (1)
and (3) we have $\nu_{3a}+\nu_{11a}+\nu_{11b}=1$. Again, in (2) we
need to consider two cases: $\chi(u^3)=\chi{(11a)}$ and
$\chi(u^3)=\chi{(11b)}$, but both cases lead us to the same system
of inequalities
$$
\textstyle \mu_1(u,\chi_5) = \textstyle\frac{1}{33} ( 2\nu_{3a} +
9 )\geq 0;\quad \mu_{11}(u,\chi_5) = \textstyle\frac{1}{33} ( -20
\nu_{3a} + 9 )\geq 0,
$$
that has no integral solutions with $\mu_1(u,\chi_5), \;
\mu_{11}(u,\chi_5) \in \Bbb Z$.

\noindent $\bullet$ Let $u$ be a unit  of order $55$. Then by (1)
and (3) we have $\nu_{5a}+\nu_{11a}+\nu_{11b}=1$. Considering two
cases when $\chi(u^5)=\chi{(11a)}$ and $\chi(u^5)=\chi{(11b)}$, we
get the same systems of inequalities
$$
\textstyle \mu_5(u,\chi_8) = \textstyle\frac{1}{55} ( 4\nu_{5a} +
40 )\geq 0;\quad \mu_{5}(u,\chi_5) = \textstyle\frac{1}{55} ( -4
\nu_{5a} + 15 )\geq 0,
$$
which also has no integral solutions such that $\mu_5(u,\chi_5),
\; \mu_5(u,\chi_8) \in \Bbb Z$.

Thus, the theorem is proved, and now the corollary follows immediately.


\Refs

\ref\no{1} \by V.~A.~Artamonov, A.~A.~Bovdi
\book Integral group rings: groups of invertible elements and classical $K$-theory
\bookinfo Algebra. Topology. Geometry
\publ Translated in J. Soviet Math.
\vol 57(2)
\pages 2931--2958
\yr 1991
\endref

\ref\no{2} \by S.~D.~Berman
\paper On the equation $x\sp m=1$ in an integral group ring
\jour Ukrain. Mat. \v Z. \vol 7 \yr 1955 \pages 253--261
\endref

\ref\no{3}
\by F.~Bleher, W.~Kimmerle
\paper On the structure of integral group rings of sporadic groups
\jour LMS J. Comput. Math.
\vol 3
\yr 2000
\pages 274--306
\endref

\ref\no{4} \by V.~Bovdi, C.~H\"ofert, W.~Kimmerle\paper  On the
first Zassenhaus conjecture for integral group rings\jour Publ.
Math. Debrecen \vol 65(3--4) \yr 2004 \pages 291--303
\endref

\ref\no {5} \by  V.~Bovdi, A.~Konovalov, R.~Rossmanith,
C.~Schneider \book  LAGUNA -- Lie AlGebras and UNits of group
Algebras, Version 3.3.3 \yr 2006 \publ http://ukrgap.exponenta.ru/laguna.htm
\endref

\ref\no{6} \by J.A.~Cohn, D.~Livingstone  \paper On the
structure of  group algebras I \jour Canad. J. Math. \vol 17 \yr
1965\pages 583--593 \endref

\ref \no{7} \by  The GAP~Group \book GAP -- Groups, Algorithms,
and Programming, Version 4.4 \yr 2006 \publ http://
www.gap-system.org
\endref

\ref\no{8} \by D.~Gorenstein \book The classification of finite
simple groups.
 Vol. 1. Groups of noncharacteristic $2$ type \publ The University Series in Mathematics. Plenum Press
\publaddr New York \yr 1983 \pages 487
\endref

\ref\no{9}
\by M.~Hertweck
\paper Partial augmentations and Brauer character values of torsion
units in group rings
\toappear
\yr 2005
\pages 26 pages
\endref

\ref\no{10} \by M.~Hertweck \paper On the torsion units of some
integral group rings \jour Algebra Colloquium  \vol 13(2) \yr 2006
\pages 328--348
\endref

\ref\no{11} \by W.~Kimmerle \book On the prime graph of the unit
group of integral group rings of finite groups \bookinfo Groups,
rings and algebras. Contemporary Mathematics \publ AMS \yr to
appear
\endref

\ref\no{12} \by  I.~S.~Luthar, I.~B.~S.Passi  \paper  Zassenhaus
conjecture for $A\sb 5$\jour Proc. Indian Acad. Sci. Math. Sci
\vol 99\yr 1989 \pages  1--5
\endref

\ref\no{13} \by I.~S.~Luthar, P.~Trama \paper Zassenhaus
conjecture for $S\sb 5$\jour Comm. Algebra\vol 19(8)\yr 1991\pages
2353--2362
\endref

\ref\no{14} \by Z.~Marciniak, J.~Ritter, S.~K.~Sehgal, A.~Weiss \paper
Torsion units in integral group rings of some metabelian groups.
II \jour J. Number Theory \vol 25 \yr 1987 \pages 340--352
\endref

\ref\no{15}
\by H.~Zassenhaus
\paper On the torsion units of finite group rings
\jour Studies in mathematics (in honor of A. Almeida Costa),
Inst. Alta Cultura, Lisbon
\pages 119--126
\endref

\endRefs
\enddocument